\mag\magstep1
\hsize=6truein\hoffset=.25truein\vsize=8.9truein
\baselineskip12pt\parskip3pt plus1pt minus1pt
\mathsurround=1.5pt
\font\smc=cmcsc10 at 10truept
\font\Hfont=cmr12 at 14.4pt
\font\Authorfont=cmr10 at 12pt
\def\Hrule{\hrule width0pt height0pt}
 \def\th#1\par{\bigskip\goodbreak \noindent{\bf#1\unskip.}\sl\enspace}
 \def\endth{\goodbreak\rm}
\def\pf#1\par{\medskip\goodbreak\noindent\it#1\unskip:\rm\enspace}
\def\endpf{\goodbreak}
\def\rk#1\par{\goodbreak\noindent\bf#1\unskip.\rm\enspace}

\def\df#1\par{\goodbreak\noindent\bf#1\unskip.\rm\enspace}
\def\enddf{\goodbreak}
\def\qedbox{\hbox{\vbox{\hrule width0.2cm height0.2pt
  \hbox to 0.2cm{\vrule height 0.2cm width 0.2pt
    \hfil\vrule height0.2cm width 0.2pt}
    \hrule width0.2cm height 0.2pt}\kern1pt}}
\def\qed{\hfill\qedbox\endpf}
\def\NoindentAfter{\everypar={\setbox0=\lastbox\everypar={}}}
\long\def\H#1\par#2\par{{\baselineskip=20pt
  \parindent=0pt\parskip=0pt\frenchspacing
  \leftskip=0pt plus .2\hsize minus .3\hsize
  \rightskip=0pt plus .2\hsize minus .3\hsize
  \def\\{\unskip\break}%
  \pretolerance=10000 \Hfont #1\unskip\break\null
   \vskip7pt\Hrule \Authorfont #2\unskip\break\null\hrule}
  \vskip48pt plus 4pt minus 4pt
  \par\NoindentAfter\rm}
\long\def\HH#1\par{\bigskip\goodbreak\noindent\bf#1\par\NoindentAfter\rm}
\long\def\HHH#1\par{\medskip\goodbreak\noindent\bf#1.\quad\rm}
\def\makeheadline{%
  \ifnum\pageno=1\vbox{\line{\vbox to 8.5pt{}\hfil}}%
   \else\vbox{\ifodd\pageno\rightheadline\else\leftheadline\fi}%
   \fi\vskip 12pt}%
\def\rightheadline{\line{\vbox to 8.5pt{}%
  \smc\hfill Generalisation del Notion de Producto Libere Amalgamate\hglue5mm
\folio}}  
\def\leftheadline{\line{\vbox to 8.5pt{}%
  \unskip\smc\unskip\folio\hglue5mm John R. Stallings\hfill}}
\nopagenumbers
\def\makefootline{\ifnum\pageno=1\centerline{\smc\folio}\else\fi}
\def\lBr{\raise.125ex\hbox{[\kern.1125ex}}
\def\rBr{\raise.125ex\hbox{\kern.1125ex]}}
\def\cite#1{{\bf\lBr#1\rBr}}
\def\em#1{{\it#1\/}}
\def\hfb{\hfill\break}\def\noi{\noindent}
\def\ind#1{\noindent\null\hskip#1\parindent}
\def\inv{^{-1}}
\def\Lsh{\mathord{\raise.3ex\hbox{\hbox to 0pt{$\supset$\hss}\raise0.92ex
\hbox{$\scriptscriptstyle<$}\phantom{e}}}}
\def\wt{\widetilde}
\def\zero{\mathord{\wedge}}
\def\comp{\mathbin{\raise1pt\hbox{$\scriptstyle\circ$}}}
\def\Star{\mathbin{\raise1pt\hbox{$\star$}}}
\H Sur un generalisation del notion de\\ producto libere amalgamate de gruppos
\par
John R. Stallings
\par
\HH Introduction
\par
In \cite{1}, io ha studite lo que occurre in le circumstantia que un gruppo $G$
ha un subensemble $P$ tal que tote elemento de $G$ es representabile unicamente
per un verbo reducite in $P$.  Il eveni que tal $G$ es multo como un producto
libere.
\par
Que occurre quando le representation per verbo reducite es unic solmente modulo
le sorta de equivalentia que interveni in le theoria del productos libere
amalgamate?  In iste articulo, io determina le structura internal del
subensemble $P$ (io los appella ``pregruppos''), e prova, sequente le methodo
de van der Waerden, que su gruppo universal ha le proprietate desiderate. 
Multe interessante exemplos pote esser trovate; tote semble simile aliquanto al
productos libere amalgamate; sed il es nulle simple maniera de construer los
omne ex ordinari tal productos.
\par
\HH1. Definition e enunciation del theorema
\par
\df1.1. Definition
\par
Un \em{pregruppo} consiste de:\hfb
\ind2$(a)$ Un ensemble $P$.\hfb
\ind2$(b)$ Un elemento de $P$, denotate per $1$.\hfb
\ind2$(c)$ Un function $P\to P$, denotate per $x\mapsto x\inv$.\hfb
\ind2$(d)$ Un subensemble $D$ de $P\times P$.\hfb
\ind2$(e)$ Un function $D\to P$, denotate per $(x,y)\mapsto xy$.\hfb
Tal que le cinque axiomes sequente sia ver:\hfb
\ind1$(1)$ Pro tote $x\in P$ es que $(1,x)$, $(x,1)\in D$ e que $1x=x1=x$.\hfb
\ind1$(2)$ Pro tote $x\in P$ es que $(x,x\inv)$, $(x\inv,x)\in D$ e que\hfb
\ind2$xx\inv=x\inv x=1$.\hfb
\ind1$(3)$ Pro tote $x,y\in P$, si $(x,y)\in D$, alora $(y\inv,x\inv)\in D$ e\hfb
\ind2$(xy)\inv=y\inv x\inv$.\hfb
\ind1$(4)$ Pro tote $x,y,z\in P$, si $(x,y)$, $(y,z)\in D$, alora:  $(x,yz)\in D$ si
e solo\hfb
\ind2si $(xy,z)\in D$, in qual caso, es que $x(yz)=(xy)z$.\hfb
\ind1$(5)$ Pro tote $w,x,y,z\in P$, si $(w,x)$, $(x,y)$, $(y,z)\in D$, alora\hfb
\ind2o $(w,xy)\in D$ o $(xy,z)\in D$.
\enddf
\par
Nos dice subinde que $xy$ es \em{definite}, in loco de $(x,y)\in D$.
\par
\df1.2. Definition
\par
Sia $P$ un pregruppo.  Un \em{verbo} in $P$ es un $n$-ena, pro alicun $n\ge 1$,
de elementos de $P$, assi:  $(x_1,\dots,x_n)$.  Le numero $n$ es appellate le
\em{longitude} del verbo.  Il es possibile \em{reducer} le verbo
$(x_1,\dots,x_n)$, si, pro alicun $i$, il es $x_i x_{i+1}$ definite; alora
$(x_1,\dots,x_{i-1},x_i x_{i+1}, x_{i+2},\dots,x_n)$ es appellate un de su
\em{reductiones}.  Le verbo es dicite esser reducite si nulle reduction existe;
i.e., pro tote $i$, es que $(x_i,x_{i+1})\notin D$.  Tote verbo de longitude un
es reducite.
\par
Si $(x_1,\dots,x_n)$, $(a_1,\dots, a_{n-1})$ es verbos, e (conveninte que
$a_0=a_n=1$) le productos $x_ia_i$, $a_{i-1}\inv x_i, a_{i-1}\inv x_i a_i$ es
toto definite, alora nos defini le \em{interfoliation} de illo per isto, a
esser:
$$(x_1,\dots,x_n)\Star (a_1,\dots,a_{n-1})=(y_1,\dots,y_n)$$
que $y=a_{i-1}\inv x_ia_i$.
\enddf
\par
Nos provera:
\par
\medskip\goodbreak\sl
$(1)$ Si $X$ es reducite e le interfoliation $X\Star A$ es definite, alora
$X\Star A$ es reducite.
\par
$(2)$ Le relation sur verbos reducite que $X\approx X\Star A$ es un relation
de equivalentia.
\par
\medskip\goodbreak\rm
Pro $a\in P$ e $X=(x_1,\dots,x_n)$ un verbo reducite, defini $\lambda_a(X)$
assi:
\par
\ind1Si $ax_1$ es non definite,
$$\lambda_a(X)=(a,x_1,\dots,x_n).$$
\ind1Si $ax_1$ es definite sed $(ax_1)x_2$ non definite,
$$\lambda_a(X)=(ax_1,x_2,\dots,x_n).$$
\ind1Si $ax_1$ e $(ax_1)x_2$ es definite,
$$\lambda_a(X)=((ax_1)x_2,x_3,\dots, x_n).$$
\par
\medskip\goodbreak\sl
\par
$(3)$ Si $X$ es reducite, etiam es $\lambda_a(X)$.
\par
$(4)$ Si $X$ es reducite e $ab$ es definite, alora $\lambda_{ab}(X)\approx
\lambda_a(\lambda_b(X))$, quo $\approx$ es definite in $(2)$.
\par
$(5)$ Le pregruppo $P$ pote esser incorporate in un gruppo universal $U(P)$,
tal que tote elemento $g\in U(P)$ pote esser scribite como un producto
$$g=x_1 x_2 \cdots x_n$$
quo $(x_1,\dots,x_n)$ es un verbo reducite in $P$, e duo tal verbos reducite
pro le mesme $g$ es $\approx$ equivalente.
\medskip\goodbreak\rm
\par
Le theorema $(5)$ es le principal; le alteres es lemmas pro isto, in le maniera
del prova per van der Waerden \cite{3} pro un resultato simile in le contexto
del productos libere.  Illo ha le corollario que \em{tote pregruppo es
continite fidelmente in su gruppo universal}.  Existe multe exemplos de
pregruppos, de que un dona le producto libere con amalgamation; un altere dona
lo que nos ha denotate per $A_F{\Lsh}\varphi$ \cite{2}.  E alteres existe plus
estranie.  Ante doner tal exemplos, nos provera le theorema.
\HH2.  Lemmas (1) (2) (3)
\par
In iste section e le proxime, sia $P$ un pregruppo fixe.
\par
\th2.1
\par
$(x\inv)\inv=x$.
\endth
\pf Prova
\par
Applica axiomes $(4)$, $(2)$, e $(1)$ al producto $xx\inv(x\inv)\inv$.\qed
\par
\th2.2
\par
Si $ax$ es definite, alora $a\inv(ax)$ es definite, e $a\inv(ax)=x$. 
Dualmente, si $xa$ es definite, etiam es $(xa)a\inv$, e $(xa)a\inv=x$.\endth
\pf Prova
\par
Par axioma $(2)$, es que $a\inv a$ es definite e $=1$.  Dunque, per
axiomas $(4)$ e $(1)$ es $a\inv(ax)$ definite e $=(a\inv a)x=x$.  Le caso dual
se prova mesmo.\qed
\par
\th2.3
\par
Si $xa$ e $a\inv y$ es definite, alora:  $xy$ es definite si e solo si
$(xa)(a\inv y)$ es definite, in qual caso $xy=(xa)(a\inv y)$.\endth
\pf Prova
\par
Applica axioma $(4)$ e 2.2 al producto de $x,a,(a\inv y)$.\qed
\par
\th2.4
\par
Si $xa$ e $a\inv y$ es definite, alora $(x,y,z)$ es un reducite verbo si e solo
si\hfb $(xa,a\inv y, z)$ es reducite.  Dualmente, $(z,x,y)$ es reducite si e solo
si $(z,xa,a\inv y)$ es reducite.\endth
\pf Prova
Nos debe monstrar que si $(x,y,z)$ es reducite, alora $(a\inv y)z$ es non
definite.  Suppone que $(a\inv y)z$ es definite e considera $\{x,a,a\inv y,z\}$. 
Alora, per axioma $(5)$ (e 2.2 e 2.1, per provar que $a(a\inv y)$ es definite),
o $x(a(a\inv y))$ es definite o $(a(a\inv y))z$ es definite.  Perque $a(a\inv
y)=y$, in ambe casos $(x,y,z)$ es non reducite.  Dunque, $(a\inv y)z$ es non
definite si $(x,y,z)$ es reducite; per 2.3, es $(xa)(a\inv y)$ anque non
definite; ita es $(xa,a\inv y, z)$ reducite.
\par
Le converso e le dual se prova mesmo.\qed
\par
On pote provar que axiomas $(1)$--$(4)$ con 2.4 implica axioma $(5)$.  Dunque,
perque axiomas $(1)$--$(4)$ es rationabile e natural, il necessita axioma $(5)$
pro nostre investigation.
\par
\th 2.5
\par
Si $(x,y)$ es un verbo reducite, e si $xa$, $a\inv y$, $yb$ es definite, alora
$(a\inv y)b$ es definite.\endth
\pf Prova
\par
Si non, per 2.3 es $(xa,a\inv y, b)$ reducite.  Dunque, per 2.4, es $(x,y,b)$
reducite, in contradiction con $yb$ esser definite.\qed
\th2.6
\par
Si $(x,y)$ es un verbo reducite, si $xa$, $a\inv y$, $(xa)b$, $b\inv(a\inv y)$
es definite, alora $ab$ es definite.\endth
\pf Prova
\par
Per 2.3 bis, $((xa)b)(b\inv(a\inv y))$ es non definite.  Applica axioma $(5)$
pro $\{x\inv,xa,b,b\inv(a\inv y)\}$; le consecutive productos es definite per
2.2; le producto del ultime trina es non definite.  Dunque per axioma $(5)$, le
producto del prima trina es definite.  Per axioma $(4)$ es
$x\inv((xa)b)=(x\inv(xa))b=ab$, per 2.2, definite.\qed
\par
\th 2.7. Lemma (1)
\par
Sia $X=(x_1,\dots,x_n)$ un verbo reducite, e $A=(a_1,\dots,a_{n-1})$ un verbo. 
Conveni que $a_0=a_n=1$.  Suppone que $x_ia_i$ e $a_{i-1}\inv x_i$ es
definite.  Alora $(a_{i-1}\inv x_i)a_i$ es definite, e $Y=X\Star A=(x_1a_1,
a_1\inv x_2 a_2,\dots,a_{n-1}\inv x_n)$ es reducite.\endth
\par
\pf Prova
\par
Applica 2.4 e 2.5 al subverbos de $X$; 2.5 monstra que $(a_{i-1}\inv x_i)a_i$
es definite.  A fortia de axioma $(4)$, nos pote omitter le parentheses.  2.4
monstra que $X\Star A$ es reducite.\qed
\par
\df2.8. Definition
\par
Per $R_n$ o $R_n(P)$, nos denota le ensemble del verbos reducite de $P$ del
longitude $n$.  Per $P^{n-1}$, nos denote le ensemble de tote verbos de $P$ del
longitude $n-1$.  Si $A=(a_1,\dots,a_{n-1})$ e $B=(b_1,\dots,b_{n-1})\in
P^{n-1}$, e si $a_i b_i$ es definite pro tote $i$, alora nos denota per $AB$ le
verbo $(a_1b_1,\dots,a_{n-1}b_{n-1})$.\enddf
\par
\th 2.9
\par
Si $X\in R^n$, $A,B\in P^{n-1}$, e si $X\Star A$ e $(X\Star A)\Star B$ pote
esser definite, alora $AB$ pote esser definite, e pois
$$(X\Star A)\Star B=X\Star(AB)$$\endth
\par
\pf Prova
\par
Applica 2.6 al subverbos de $X$.  Isto monstra que $AB$ pote esser definite. 
Axiomas $(4)$ e $(3)$ monstra que $(X\Star A)\Star B=X\Star(AB)$.\qed
\par
\df2.10. Definition
\par
Le relation $\approx$ sur $R_n$ es definite assi:
$$(x_1,\dots,x_n)\approx(y_1,\dots,y_n)$$
si e solo si existe $(a_1,\dots,a_{n-1})\in P^{n-1}$ tal que cata $x_ia_i$ e
$a_{i-1}\inv x_i$ es definite, e $y_i=a_{i-1}\inv x_i a_i$.  I.e., $X\approx Y$
si e solo si existe $A$ tal que $Y=X\Star A$.\enddf
\th2.11.  Lemma (2)
\par
Le relation $\approx$ es un relation de equivalentia sur $R_n$.\endth
\par
\pf Prova
\par
Si $I=(1,\dots,1)$, alora $X= X\Star I$, ita $X\approx X$.  Si
$A=(a_1,\dots,a_{n-1})$ sia $A\inv=(a_1\inv,\dots,a_{n-1}\inv)$; alora,
$Y=X\Star A$ si e solo si $X=Y\Star A\inv$; dunque, si $X\approx Y$, alora
$Y\approx X$.  Si $Y=X\Star A$ e $Z=Y\Star B$, alora per 2.9, es $AB$
definabile e $Z=X\Star (AB)$; dunque, si $X\approx Y$ e $Y\approx Z$, alora
$X\approx Z$.\qed
\par
\df2.12. Definition
\par
Per $R$ o $R(P)$, nos denota le reunion de tote $R_n$, pro $n=1,2,3,\dots$.  Pro
cata $a\in P$ e cata $X\in R$, nos defini un verbo $\lambda_a(X)$ como seque:
\par
Sia $X=(x_1,x_2,x_3,\dots)$.\hfb
\ind1$(1)$ Si $(a,x_1)$ es reducite, alora
$$\lambda_a(X)=(a,x_1,x_2,\dots).$$
\ind1$(2)$ Si $ax_1$ es definite sed $(ax_1,x_2)$ reducite, alora
$$\lambda_a(X)=(ax_1,x_2,x_3,\dots).$$
\ind1$(3)$ Si $ax_1$ e $(ax_1)x_2$ es definite, alora
$$\lambda_a(X)=((ax_1)x_2,x_3,\dots).$$
In caso $(2)$ nos includa, como caso degenerate, le possibilitate que $X$ ha
longituda un, quando $ax_1$ es definite.\enddf
\th2.13. Lemma (3)
\par
Si $X$ es reducite, alora $\lambda_a(X)$ es reducite.\endth
\pf Prova
\par
Isto es obvie in le casos $(1)$ e $(2)$.  In caso $(3)$, quo $ax_1$,
$(ax_1)x_2$ es definite, sed $x_1x_2$ e $x_2 x_3$ non definite, nos debe provar
que $((ax_1)x_2)x_3$ es non definite.  Considera $\{x_1,x_1\inv
a\inv,(ax_1)x_2,x_3\}$ e lo applica axioma $(5)$; isto es possibile si
$((ax_1)x_2)x_3$ es definite; mais axioma $(5)$ implicarea alora que o $x_1x_2$
o $x_2x_3$ es definite, in contradiction con $X$ esser reducite.\qed
\par
\HH3. Lemma (4)
\par
Hic nos prova
\th 3.1.  Lemma (4)
\par
Si $X$ es reducite e $ab$ es definite, alora $\lambda_{ab}(X)\approx
\lambda_a(\lambda_b(X))$.\endth
\pf Prova
\par
Le prova consiste de spectar le varie casos.  Sia $X=(x_1,\dots,x_n)$.
\par
\em{Caso} $1$:  Es $bx_1$ non definite.  Alora,
$$\lambda_b(X)=(b,x_1,\dots,x_n).$$
\ind1\em{Subcaso} $1_1$:  Es $(ab)x_1$ non definite.  Pro applicar $\lambda_a$
nos nos trove in caso $2.12(2)$; dunque
$$\lambda_a(\lambda_b(X))=(ab,x_1,\dots,x_n)=\lambda_{ab}(X).$$
\ind1\em{Subcaso} $1_2$:  Es $(ab)x_1$ definite.  Alora, pro applicar
$\lambda_a$ a $\lambda_b$ nos nos trova in caso $2.12(3)$; dunque
$$\lambda_a(\lambda_b(X))=((ab)x_1,x_2,\dots,x_n).$$
Il seque que $((ab)x_1)x_2$ es non definite; ita, pro applicar $\lambda_{ab}$ a
$X$ nos nos trova in caso $2.12(2)$; dunque
$$\lambda_{ab}(X)=((ab)x_1,x_2,\dots,x_n)=\lambda_a(\lambda_b(X)).$$
\ind1\em{Caso} $2$:  Es $bx_1$ definite sed $(bx_1)x_2$ non definite.  Alora,
$$\lambda_b(X)=(bx_1,x_2,\dots,x_n).$$
\ind1\em{Subcaso} $2_1$:  Es $a(bx_1)$ non definite.  Alora, $(a,bx_1)$ es
reducite; dunque\hfb $(ab,b\inv(bx_1))=(ab,x_1)$ es redu;cite.  In iste caso,
$$\eqalign{\lambda_a(\lambda_b(X))&=(a,bx_1,x_2,\dots,x_n),\cr
\lambda_{ab}(X)&=(ab,x_1,x_2,\dots,x_n).\cr}$$
Dunque, $\lambda_{ab}(X)=(\lambda_a(\lambda_b(X)))\Star(b,1,\dots,1)$.\hfb
\ind1\em{Subcaso} $2_2$:  Es $(a(bx_1)$ definite.  Alora, anque $(ab)x_1$ es
definite e $=a(bx_1)$.  Si $(abx_1)x_2$ es non definite, alora
$$\lambda_{ab}(X)=\lambda_a(\lambda_b(X))=(abx_1,x_2,\dots,x_n).$$
Si $(abx_1)x_2$ es definite, alora
$$\lambda_{ab}(X)=\lambda_a(\lambda_b(X))=((abx_1)x_2,x_3,\dots,x_n).$$
\ind1\em{Caso} $3$:  Es $bx_1$ e $(bx_1)x_2$ definite.  Alora,
$$\lambda_b(X)=((bx_1)x_2,x_3,\dots,x_n).$$
\ind1\em{Subcaso} $3_1$:  Es $a(bx_1)$ non definite.  Alora, $(ab)x_1$ es non
definite, e
$$\lambda_{ab}(x_1,\dots,x_n)=(ab,x_1,\dots,x_n)\approx
(a,bx_1,x_2,\dots,x_n).$$
Ita, $(bx_1)x_2$ es non definite.  In altere verbos, iste subcaso non occurre
jammais.\hfb
\ind1\em{Subcaso} $3_2$:  Es $a(bx_1)$ definite.  Alora, $(ab)x_1$ es definite e
$=a(bx_1)$.\hfb
\ind1\em{Subsubcaso} $3_{2_1}$:  Es $(abx_1)x_2$ non definite.  Alora
$a((bx_1)x_2)$ es non definite.  Dunque:
$$\eqalign{\lambda_{ab}(X)&=(abx_1,x_2,\dots,x_n),\cr
\lambda_a(\lambda_b(X))&=(a,(bx_1)x_2,x_3,\dots,x_n).\cr}$$
Ita
$$\lambda_{ab}(X)=(\lambda_a(\lambda_b(X)))\Star (bx_1,1,\dots,1).$$
\ind1\em{Subssubcaso} $3_{2_2}$:  Es $(abx_1)x_2$ definite.  Alora
$(a(bx_1))x_2=a((bx_1)x_2)$ es definite, e
$$\lambda_{ab}(X)=((abx_1)x_2,x_3,\dots,x_n).$$
Dunque, $((abx_1)x_2)x_3$ es non definite; ita,
$$\lambda_a(\lambda_b(X))=((abx_1)x_2,x_3,\dots, x_n)=\lambda_{ab}(X).$$
\par
Isto exhauri tote le casos possibiles, ita prova $3.1$.\qed
\par
\HH 4. Le Theorema Principal
\par
Pro un pregruppo $P$, nos ha $R(P)$, le ensemble de tote reducite verbos in
$P$, sur le qual nos ha le relation de equivalentia $\approx$.  Per $\wt R(P)$,
nos denota le ensemble del classes de $\approx$ equivalentia.
\par
\th4.1
\par
Pro tote $a\in P$, le function $\lambda_a:R(P)\to R(P)$ induce un function,
anque denotate $\lambda_a$, $\wt R(P)\to\wt R(P)$.\endth
\pf Prova
\par
Nos debe monstrar que, si $Y=X\Star B$, alora $\lambda_a(X)\approx
\lambda_a(Y)$.  Sia $X=(x_1,\dots,x_n)$ e $B=(b_1\dots,b_{n-1})$.  Le tres
casos es como seque:\hfb
\ind1$(1)$ $ax_1$ es non definite.  Sia $B'=(1,b_1,\dots,b_{n-1})$.  In iste
caso, $\lambda_a(Y)=(\lambda_a(X))\Star B'$.\hfb
\ind1$(2)$ $ax_1$ es definite sed $(ax_1)x_2$ non definite.  In iste caso,
$\lambda_a(Y)=(\lambda_a(X))\Star B$.\hfb
\ind1$(3)$ $ax_1$ e $(ax_1)x_2$ es ambes definite.  Sia
$B''=(b_2,\dots,b_{n-1})$.  In iste caso, $\lambda_a(Y)=(\lambda_a(X))\Star
B''$.
\par
In cata caso, le verbo $\lambda_a(X)$ es determinate; le expression al latere
dextre del formula determina un verbo reducite; examine de iste verbo determina
$\lambda_a(Y)$.  Per exemplo, in caso $(3)$,
$$\eqalign{\lambda_a(X)&=((ax_1)x_2,x_3,\dots,x_n)\cr
\lambda_a(X)\Star B''&=((ax_1)x_2b_2,b_2\inv x_3b_3,\dots,b_{n-1}\inv x_n)\cr
&=((ax_1b_1)b_1\inv x_2b_2,\dots,b_{n-1}x_n)\cr
&=((ay_1)y_2,\dots,y_n)\cr
&=\lambda_a(Y).\cr}$$
Le producto triple $(ax_1)x_2b_2$ es definite, per le dual de $2.5$, perque
$(x_2,x_3)$ es reducite, $x_2b_2$ e $b_2\inv x_3$ es definite, e $(ax_1)x_2$ es
definite.  Mesmo es $ax_1b_1$ definite.
\par
Le detalios del altere casos es plus facile.\qed
\df4.2. Definition
\par
Sia $P$, $Q$ duo pregruppos.  Un function $\varphi:P\to Q$ es appellate un
\em{morphismo} de pregruppos si, pro tote $x,y\in P$ tal que $xy$ es definite,
es que $\varphi(x)\varphi(y)$ es definite, e que
$\varphi(xy)=\varphi(x)\varphi(y)$.
\par
La classe de pregruppos con lor morphismos constitue un categoria continente le
categoria de gruppos e homomorphismos.  Il existe, per nonsenso abstracte, un
functor co-adjuncte al functor de inclusion.  Isto dona le \em{gruppo
universal} $U(P)$ de un pregruppo.  I.e.,  $U(P)$ es un gruppo, e ha un
morphismo specific $\iota:P\to U(P)$, tal que, pro tote gruppo $G$ e tote
morphismo $\varphi:P\to G$, il existe un homomorphismo unic $\psi:U(P)\to G$ tal
que $\varphi=\psi\comp \iota$.\enddf
\th4.3
\par
Sia $S$ le gruppo de permutationes de $\wt R(P)$.  Alora $\lambda$ es un
morphismo de $P$ in $S$.\endth
\pf Prova
\par
Perque $\lambda_1$ es le function identic de $\wt R(P)$ in se, e que
$\lambda_x\comp\lambda_{x\inv}=\lambda_{x\inv}\comp \lambda_x=\lambda_1$, per
$3.1$, tote $\lambda_x$ pertine a $S$.  Anque, $3.1$ pote esser interpretate
como dicer $\lambda$ esser un morphismo.\qed
\noi\bf4.4.\rm\enspace Per le proprietate universal, $\lambda$ extende
unicamente a un homomorphismo, anque denotate $\lambda$, de $U(P)$ in $S$.  Nos
denota le valor de $\lambda$ sur $g\in U(P)$, per $\lambda_g:\wt R(P)\to\wt
R(P)$.
\par
Perque $\iota(P)$ genera $U(P)$, cata $g\in U(P)$ pote esser scribite como
$g=\iota(x_1)\iota(x_2)\cdots\iota(x_n)$, quo $(x_1,x_2,\dots,x_n)$ es un verbo
in $P$.  Post applicar reductiones a iste verbo, nos obtene un verbo reducite
$(x_1,\dots,x_n)$ tal que $g=\iota(x_1)\iota(x_2)\cdots \iota(x_n)$.  Nos
denota per $\zero$, le verbo $(1)$ de longitude un.  Nos ha le formula:
$$\lambda_g([\zero])=%
\lambda_{x_1}(\lambda_{x_2}(\cdots(\lambda_{x_n}([\zero]))\cdots))=%
[(x_1,x_2,\dots,x_n)],$$
quo $[\phantom{\zero}]$ denota le classe de $\approx$ equivalentia.  Tote
application de un $\lambda_{x_i}$ hic se trova in le caso $2.12(1)$.
\par
In iste maniera $g$ determina per se le classe del verbos reducite que
representa $g$.  Dunque:
\th4.5. Theorema
\par
Si $P$ es un pregruppo, alora tote elemento $g\in U(P)$, le gruppo universal de
$P$, pote esser representate como producto $x_1x_2\cdots x_n$ de un verbo
reducite in $P$, $(x_1,\dots,x_n)$.  Duo tel verbos reducite representa le
mesme elemento de $U(P)$, si et solo si illos es $\approx$ equivalente.  \rm(Hic
nos ha identificate $x\in P$ con $\iota(x)\in U(P)$, pro simplificar le
notation.)\endth
\pf Prova
\par
Nos ha jam provate le ``solo si''.  Le ``si'' es le computation trivial que, si
$g=x_1x_2\cdots x_n$, alora  $g=(x_1a_1)(a_1\inv x_2a_2)\cdots(a_{n-1}\inv x_n)$.
\qed
\th4.6. Corollario
\par
Un pregruppo $P$ es continite fidelmente in su gruppo universal $U(P)$.
\endth
\pf Prova
\par
Isto vole dicer que le morphismo specific $\iota:P\to U(P)$ es injective.  Isto
seque del theorema, perque nulle verbo de longitude un non es equivalente a
nulle altere verbo.\qed
\par
\HH5. Exemplos
\par
\df5.1
\par
Le plus standard exemplo de un pregruppo es facite de tres gruppos $A$, $B$,
$C$ e de duo monomorphismos $\varphi:C\to A$, $\psi:C\to B$.  Identifica
$\varphi(C)$ con $\psi(C)$; alora $A\cap B=C$.  Sia $P=A\cup B$.  Le $1$ e le
inverso es obvie; le producto es definite pro duo elementos $x,y$ si e solo si
le duo pertine a un singule del $A$ o $B$.  Le axiomas $(1)$ usque $(4)$ es
clarmente satisfacte.  Pro axioma $(5)$, il frange in casos simple facile a
verificar.  Le gruppo universal es le producto libere con amalgamation $A*_CB$.
\enddf
\df5.2
\par
Ecce un caso simile sed plus general.  Un \em{arbore de gruppos} consista
de:\hfb
\ind1$(a)$ Un ensemble $I$, partialmente ordinate per $<$, con elemento minime,
tal que pro tote $i,j,k\in I$, si $i<k$ e $j<k$, alora o $i\le j$ o $j\le i$. 
(Tal ensemble ordinate es un sorta de arbore abstracte.)\hfb
\ind1$(b)$ Un classe de gruppos $\{G_i\}$ indicate per $i\in I$.\hfb
\ind1$(c)$  Per tote $i,j\in I$, si $i<j$, un monomorphismo $\phi_{ij}:G_i\to
G_j$; tal que, per tote $i,j,k\in I$, si $i<j<k$, alora $\phi_{jk}\comp
\phi_{ij}=\phi_{ik}:G_i\to G_k$.
\par
Nos pote construer, como supra, le reunion $P$ de tote $\{G_i\}$, identificante
$x\in G_i$ con $\phi_{ij}(x)\in G_j$.  Le lector pote verifica que, a fortia
del proprietates de arbore, con le obvie operationes, $P$ es un pregruppo.  Le
gruppos universal de tal progruppos include tote ordinari productos libere con
amalgamation de multe factores.\enddf
\df5.3
\par
Considera un producto libere amalgamate $A*_CB$.  Sia $P$ le subensemble de
tote elementos que pote esser scribite $bab'$, pro alicun $b,b'\in B$, $a\in
A$; dunque, $P$ contine $A$ e $B$ e aliquanto plus.  Dice que la producto $xy$
de duo elementos $x,y\in P$ es definite, quando $xy\in P$.  Usante le structura
(per verbos reducite in $A\cup B$, etc.) de $A*_CB$, nos pote provar que $P$
es un pregruppo.  Le gruppo universal de $P$ es etiam $A*_CB$; mais le
structura de $A*_CB$ per verbos in $P$ es differente de illo per verbos in
$A\cup B$.\enddf
\df5.4
\par
Considera un gruppo $G$ con subgruppo $H$.  Sia $P$ le ensemble $G$, sed defini
multiplication de $x$ e $y$ si e solo si al minus un de $\{x,y,xy\}$ pertine a
$H$.  Isto es un pregruppo, e su gruppo universal es non troppo simile a un
producto libere amalgamate.\enddf
\df5.5
\par
Sia $G$ un gruppo, $H$ un subgruppo, e $\varphi:H\to G$ un monomorphismo. 
Construe quatro ensembles,
$$G,\ x\inv G,\ Gx,\ x\inv Gx.$$
Identifica $h\in H\subset G$, con $x\inv\varphi(h)x\in x\inv Gx$.  Defini
multiplication inter $G$ e $G$, $G$ e $Gx$, $x\inv G$ e $G$, $x\inv G$ e $Gx$,
$Gx$ e $x\inv G$, $Gx$ e $x\inv Gx$, $x\inv Gx$ e $x\inv G$, $x\inv Gx$ e
$x\inv Gx$, per cancellation de $xx\inv$ e multiplication in $G$.  Per le
formulas:
$$\eqalign{hx\inv&=x\inv\varphi(h)\cr
xh&=\varphi(h)x,\cr}$$
que seque del identification de $H$ con $x\inv\varphi(H)x$, multiplication es
defini in tote caso quando un factor pertina a $H$.  Iste monstruositate es un
pregruppo, le gruppo universal de que es appellate  $G_H\Lsh \varphi$.\enddf
\vfill\eject
\noi\centerline{\bf Referentias}
\medskip
\noi\cite{1}. J. Stallings, ``A remark about the description of free products
of groups'', Proc.~Cambridge Philos.~Soc.~{\bf62} (1966), 129--134.
\par
\noi\cite{2}. J. Stallings, ``On the theory of ends of groups'', (a parer).
\par
\noi\cite{3}.  B.L. van der Waerden.``Free products of groups'',
Amer.~J.~Math.~{\bf70} (1948), 527--528.
\vskip.1in
\noi\line{*\hfil*\hfil*\hfil*\hfil*\hfil*\hfil*\hfil*\hfil*\hfil*\hfil*\hfil*}
\par
\noi\line{\hfill Universitate de California}
\line{\hfill Berkeley}
\line{\hfill februario 1968}
\vskip.1in
\hrule
\vskip.1in
Vinti-cinque annos retro, io ha scribite iste articulo.  Es in le lingua 
international ``Interlingua'', como describite in le libro del IALA ({\sl
Interlingua\/}, 1951, Storm Publishers, New York), que io ha emite in 
Telegraph Avenue in le anno ante.  Interlingua es descendite de ``Latino 
sine Flexione'' que era usate per G. Peano in su scripturas mathematic.
\par
Sub le arbores del Ca\~non del Fragas, io habeva contemplate le problema de 
Serre si un gruppo sin torsion que contine un subgruppo libere de indice 
finite esserea libere.  Isto fructava in mi opera sur le fines del gruppos
({\sl Group Theory and Three-Dimensional Manifolds\/}, 1971, Yale University
Press).  Le  notion de pregruppo era parte de isto.  Era presagite in opera
de R. Baer sur le lege associative (Amer.~J.~Math., 1949--50).  Le duo
exemplos principal,  $A*_CB$ e  $G_H\Lsh \varphi$, es debite,
respectivemente, a Schreier (Hamburg.~Abh., 1927) e a Higman-Neumann-Neumann
(J.~London Math.~Soc., 1949).
\par
\vskip.1in
\noi\line{*\hfil*\hfil*\hfil*\hfil*\hfil*\hfil*\hfil*\hfil*\hfil*\hfil*\hfil*}
\par
\noi\line{John R. Stallings\hfill Universitate de California}
\line{e-posta: {\tt stall@math.berkeley.edu}\hfill Berkeley}
\line{\hfill junio 1993}
\bye